%&amstex

\input amstex.tex
\documentstyle{amsppt}

\topmatter
\title Cores of s-cobordisms of 4-manifolds 
\endtitle
\author Frank Quinn\endauthor
\thanks Partially supported by the National Science Foundation\endthanks
\date March 2004\enddate
\address Mathematics, Virginia Tech, Blacksburg VA 24061-0123\endaddress
\email quinn\@math.vt.edu\endemail
\abstract The main result is that an s-cobordism (topological or smooth) of 4-manifolds has a product structure outside a ``core'' sub s-cobordism.  These cores are arranged to have quite a bit of structure, for example they are smooth and abstractly (forgetting boundary structure) diffeomorphic to a standard neighborhood of a 1-complex. The decomposition is highly nonunique so cannot be used to define an invariant, but it shows the topological s-cobordism question reduces to the core case.  The simply-connected version of the decomposition (with 1-complex a point) is due to Curtis, Freedman, Hsiang and Stong. Controlled surgery is used to reduce topological triviality of core s-cobordisms to a question about controlled homotopy equivalence of 4-manifolds. There are speculations about further reductions. \endabstract
\endtopmatter

\document
\head 1. Introduction\endhead
 The classical s-cobordism theorem asserts that an s-cobordism of $n$-manifolds (the bordism itself has dimension $n+1$) is isomorphic to a product if $n\geq 5$. ``Isomorphic'' means smooth, PL or topological, depending on the structure of the s-cobordism. In dimension 4 it is known that there are smooth s-cobordisms without smooth product structures; existence was demonstrated by Donaldson \cite{3}, and  specific examples identified by Akbulut \cite{1}. In the topological case product structures follow from disk embedding theorems. The best current results require ``small'' fundamental group, Freedman-Teichner \cite{5},  Krushkal-Quinn \cite{9} so s-cobordisms with these groups are topologically products. The  large fundamental group question is still open.
 
Freedman has developed several link questions equivalent to the 4-dimensional ``surgery conjecture'' for arbitrary fundamental groups. However ``surgery'' is equivalent to disk embeddings in which the manifold is allowed to change up to s-cobordism, so these link formulations do not offer insights into s-cobordisms. The objective of this paper is to begin development of primitive  questions that could detect failure of s-cobordisms to be topological products.

Curtis-Freedman-Hsiang-Stong \cite{2} (see also Kirby \cite{8}, Matveyev \cite{10}) showed  smooth s-cobordisms of simply-connected  4-manifolds can be given a product structure outside a contractible sub s-cobordism. In Section 2 this is extended to show a general compact s-cobordism (smooth or topological) can be given a product structure outside a ``core'' with the homotopy type of a 1-complex.  The 1-complex can be specified in advance, provided only that it map onto the fundamental group of the manifold. The proof is an application of the dual-decomposition results of Quinn \cite{13}.

Section 3 uses elementary arguments to develop properties of core s-cobordisms, illustrating the utility of the somewhat exotic definition. One consequence, following an argument of Matveyev \cite{7}, is that a product s-cobordism can be given a decomposition with arbitrarily prescribed core. More generally a core decomposition can be modified to make the core ``arbitrarily worse''.
 In particular this means the core does not provide a obstruction: global triviality does not depend on the core directly, but whether it can ``untwist'' inside the ambient manifold. This is presumably related to the fact that cores are 1-dimensional, but smooth invariants that detect nontrivial s-cobordisms depend on 2-dimensional homology classes. In the topological case the full conjecture does reduce to the conjecture for cores.
 
The fourth section collects consequences of high-dimensional s-cobordism  and surgery for core s-cobordisms. One consequence is a formulation of the topological s-cobordism question in terms of $\epsilon$ control of maps of 4-manifolds homotopy equivalent to graphs. We suggest ways this might lead to a further reduction of the problem. 

\head 2. Cores in s-cobordisms\endhead
Core s-cobordisms are defined in 2.1. The main theorem is given as 2.2, and the proof occupies the rest of the section.

\subhead 2.1 Definition\endsubhead 
A {\it core s-cobordism\/} consists of:
\roster\item $W^5$, a smooth 5-dimensional regular neighborhood of a graph;
\item a decomposition $\partial W^5 = N_0\cup N_1$, where $N_i$ are smooth submanifolds  intersecting in their common boundary, denoted $\partial N$; and
\item the $N_i$ have handlebody decompositions with only 0- 1- and 2-handles and spines that homotopy 2-deform to 1-complexes.
\endroster
Explanations and consequences of this definition, particularly the spine condition in (3), are given in Section 3.

We consider a core $(W^5,N_0,N_1)$ as an s-cobordism from $N_0$ to $N_1$ with a product structure given on the boundary.  The product boundary comes from inserting a collar between the pieces:
$\partial W^5 \simeq N_0\cup_{\partial N\times\{0\}} \partial N\, \times I\cup_{\partial N\times\{1\}}N_1$. 

\proclaim{2.2 Theorem} Suppose $V$ is a compact s-cobordism of 4-manifolds $M_0$, $M_1$, with a product structure on the s-cobordism of boundaries, and suppose $Y\subset M_0$ is an embedded 1-complex surjective on $\pi_1$. Then there is a decomposition $M_0=U\cup N_0$ with $Y\subset N_0$ so that:
\roster\item there is a core s-cobordism $(W^5,N_0,N_1)$ and an extension of the decomposition of $M_0$ to an isomorphism $U\times I\cup_{\partial N\times I}W^5\simeq V$ that extends the given product structure over $\partial M$;
\item the inclusions $Y\subset N_0$ and $U\subset M_0-Y$ are homotopy equivalences.
\endroster\endproclaim
``Isomorphism'' here means either diffeomorphism or homeomorphism, depending on the structure of the initial s-cobordism, Note the core is smooth whether the s-cobordism is or not. The proof of 2.2 occupies the rest of the section. Section 2.3 describes the situation halfway through the proof, where techniques change, and the end of the proof. The intermediate data is obtained in 4.4 by following the high-dimensional proof  with a little extra 4-dimensional care. In 4.5 the dual decomposition results of \cite{13} are used to complete the proof. A brief explanation of the method is given at the beginning of the section. 
\subsubhead 2.3 Intermediate data and the end of the proof\endsubsubhead
There is a decomposition $M_0=\hat U\cup \hat N_0$ that extends to a decomposition of the s-cobordism as the union of a product $U\times I$ and a cobordism $(\hat W,\hat N_0,\hat N_1)$, and  there are immersed 2-handles $H_i\:(D^2,S^1)\to (\hat U,\partial \hat U)$ so that
\roster\item Adding product handles to $\hat W$ along $(\partial H_*)\times I$ gives a product cobordism diffeomorphic to $W^4\times I$;
\item there are 2-handlebody structures on $\hat N_i$ so that adding handles on the boundaries of the $H_*$ give handlebodies that 2-deform to neighborhoods $W^4$ of $Y$;
\item the inclusion of the complement of the immersed handles $\hat U-\cup_iH_i\subset M_0$ is an isomorphism on $\pi_1$. 
\endroster
If we could find embedded disks in $\hat U$ with the same boundaries as $H_*$ then the original s-cobordism would be trivial. In Section 2.5 we show there are $\pi_1$ negligible embedded disks $\hat H_*$ in $\hat U$ with boundaries {\it homotopic\/} to the boundaries of $H_*$. Since the cobordism is  a product over $\hat U$ there are corresponding product handles in the cobordism. Add these to $\hat W$ and denote the resulting sub-cobordism by $(W,N_0,N_1)$. We now show that this gives a core decomposition.

First forgetting the  boundary structure shows the 5-manifold $W$ is obtained by adding 2-handles to circles homotopic to the attaching maps of the $H_*$. But homotopy imples isotopy for circles in a 4-manifold, so the resulting manifolds are diffeomorphic. Thus by assumption 2.3(1) $W$ is a regular neighborhood of $Y$ after deforming it into the interior. 

Next, 2.3(2) asserts that adding handles $H_*$ to the top and bottom give handlebodies that deform through 2-handlebodies to regular neighborhoods of graphs. Attaching handles on homotopic circles give manifolds with  spines the same up to homotopy of 2-cells so the spines of $N_0$ and $N_1$  homotopically 2-deform to $Y$. 

Finally let $U$ denote the complement of the interior of $N_0$ in $M_0$. By construction the complement of the interior of $W\subset V$ is a product. Since $Y\subset N_0$ is a homotopy equivalence duality shows that $M_0-Y\supset U$ is a $Z[\pi_1M_0]$ homology equivalence. But the handles deleted to give $U$ are $\pi_1$ negligible so this inclusion is also an isomorphism on fundamental groups. This implies it is a homotopy equivalence. This is the final condition needed for a core decomposition.

\subsubhead 2.4 Getting the data\endsubsubhead
The beginning of the proof of 2.2 is the same as the higher-dimensional version, with a few 4-dimensional refinements, up to the point Whitney disks are usually used. Specifically 
the proof of Theorem 7.1D in Freedman-Quinn \cite4, gives 
\roster\item  $V$ can be obtained by adding 2- and 3-handles to a collar $M_0\times I$; 
\item let $M_{1/2}$ denote the level between the handles. This has two families of  framed embedded 2-spheres: $A_i$ the attaching maps of 3-handles and $B_i$ dual spheres to the 2-handles;
\item for each $i$ there is a distinguished intersection in $A_i\cap B_i$; all other intersections between the families are arranged in pairs with immersed Whitney disks. We denote the Whitney disks by $C_*$.
\item these Whitney disks all have framed immersed transverse spheres, and the spheres and interiors of the Whitney disks are disjoint from $\cup_i(A_i\cup B_i)$; and
\item the union of $A_*$, $B_*$ and $C_*$ is $\pi_1$-negligible in the sense that removing it does not change the fundamental group of $M_{1/2}$.
\endroster
We may assume all this is disjoint from the 1-complex $Y$. The next step is to enlarge this data to include $Y$. Specifically for each $i$  choose an arc from $Y$ to the distinguished intersection in $A_i\cap B_i$. Let $\hat Y$ denote this enlarged 1-complex. We now describe how to choose an {\it accessory\/} disk similar to \cite{4, \S3.1} for each non-distinguished $A\cap B$ intersection. 

The boundary circle of an accessory disk consists of two arcs joining distinguished intersections: one in $A\cup B$ passing through the non-distinguished point, and the other in $\hat Y$. The arcs in $A\cup B$ should be embedded and disjoint from each other, except at the distinguished points, and disjoint from the Whitney arcs, except at the non-distinguished points.  Since $\pi_1\hat Y\to \pi_1M_0$ is onto the arc in $\hat Y$ can be chosen so that the resulting loops are contractible. 

Next extend these circles to 2-disks mapping to $M_0$. By the $\pi_1$ negligible hypothesis we can suppose the interiors are disjoint from $\hat Y\cup A_*\cup B_*$.  Standardize these maps near the boundary: near $A\cup B$ a collar on the circle maps to a section of the normal bundles, and near the 1-complex all these collars form a standard ``flange'' pattern. Let $\hat N_{1/2}$ be a thickening of $\hat Y\cup A_*\cup B_*$, then the standardization implies the accessory disks meet the boundary in disjoint framed embedded circles. If the framing on the boundary extends to a framing of the restriction of the tangent bundle of $M_{1/2}$ to the disk then the disk can be changed to be a framed immersion of a handle. If the framing does not extend then the disk can be ``spun'' (introducing new intersections among disks, c.f. the framing of Whitney disks in \cite{4}) to correct this, so in any case we can get represent the accessory disks as immersed 2-handles in the complement of $\text{int}(\hat N_{1/2})$ with boundaries on $\partial \hat N_{1/2}$. They can also easily be made $\pi_1$ negligible.

Whitney disks are already standardized in a neighborhood of $A\cup B$ in the sense above, so the Whitney disks also determine immersed 2-handles in the complement of $\text{int}(\hat N_{1/2})$ with boundaries on $\partial \hat N_{1/2}$. Denote the collection of Whitney and accessory disks by $H_*$. 

We now describe the cobordism $\hat W$ of the intermediate data. Recall the 3-handles of the original cobordism $V$ are attached on $A_*\subset \hat N_{1/2}$. This gives a cobordism of $ \hat N_{1/2}$ to $ \hat N_{1}\subset M_1$. Similarly the duals of the 2-handles are 3-handles attached on the spheres $B_*\subset \hat N_{1/2}$, so these handles give a cobordism of $ \hat N_{1/2}$ to $\hat N_{0}\subset M_0$. The union of these two cobordisms gives the cobordism $\hat W$. Since it contains the handles of $V$, $V$ is a product outside  $ \hat W$. 

Denote by $\hat U$ the complement in $M_0$ of the interior of $\hat N_0$. Since $V$ is a product outside $ \hat W$ we can think of the immersed handles $H_*$ in $M_{1/2} -\hat N_{1/2}$ as lying in any level, in particular in $\hat U$. 

The final step is to check what happens if we add product handles to $\hat W$ along the boundaries of the $H_*$. First note that the resulting new cobordism has embedded Whitney disks for the non-distinguished $AB$ intersection points. Pushing across these disks eliminates all these intersections. This gives 2- and 3-handles arranged in cancelling pairs, so they can be cancelled. This means the new cobordism is a product. 

We refine this picture to see the top and bottom of the cobordism. Consider a neighborhood of the union of $A_*\cup B_*$ and the embedded Whitney disks. After the Whitney moves we see this as a neighborhood of disjointly embedded dual pairs $\tilde A_i\cup \tilde B_i$, with a 1-handle attached where each Whitney move took place. Next we have the neighborhood $W^4$ of $Y$, joined by a single 1-handle to each dual pair. Finally we have the accessory disks. Each of these is attached on a circle that goes over exactly one Whitney 1-handle exactly once. Attaching 2-handles on these circles gives handles that cancel the Whitney 1-handles. The result is $W^4$ joined by 1-handles to neighborhoods of  embedded dual pairs.  This is a handle picture in the middle level $M_{1/2}$: note that it was obtained by a deformation of 2-handlebodies (i.e\. without ever using 3-handles). 

The cobordisms from the middle to $\hat N_0$ and $\hat N_1$ are obtained by attaching handles to $A_*$ and $B_*$ respectively. After adding product handles these correspond to attaching handles on half of embedded dual pairs. This  converts the dual pairs to balls. The top and bottom of the new cobordism are therefore obtained from a union of $W^4$ and some balls (0-handles) by attaching each 0-handle to $W^4$ by a single 1-handle.  These 0- and 1-handles cancel to give diffeomorphisms of the top and bottom to $W^4$. Again these reductions were accomplished by 2-handlebody moves, as required for 2.3(2). 

This completes the verification that the decomposition constructed has the properties claimed in 2.3.

\subsubhead 2.5 Homotopy handles\endsubsubhead
According to 2.3 there is a compact manifold $\hat U$ and $\pi_1$-negligible immersed 2-handles $H_*\:D^2\to \hat U$ so that if there are handles in $\hat U$ with the {\it same\/} boundaries the s-cobordism is a product, and if there are handles with {\it homotopic\/} boundaries the s-cobordism has a core. 
We find homotopic handles using the dual decompositions of \cite{13}. Since these decompositions are for smooth manifolds the first step is to reduce to the smooth case.

If $\hat U$ is not smoothable then we find a submanifold $\hat U_0\subset \hat U$ such that:
\roster\item $\hat U_0$ is smooth;
\item  $\partial \hat U\cup H_*\subset \hat U_0$; and
\item $\hat U_0\subset \hat U$ induces an isomorphism on $\pi_1$.
\endroster
The union of $\hat W$ and the product over $\hat U_0$ gives a smooth sub-s-cobordism of the original with all the same data. Therefore a solution in the smooth case applied to this sub cobordism gives a solution in the topological case. Finding such $\hat U_0$ is standard: first there is a smooth structure on the complement of a point in $\hat U$, \cite{12}. Then there is then a splitting along a homology sphere into a ``weak collar'' of the boundary and $\hat U_0$, see \cite{4, 11.9B}.

Now suppose $\hat U$ is smooth. The maps $H_*\:D^2\to \hat U$ give a chain map on cellular  chain complexes with $Z[\pi_1\hat U]$ coefficients $C_2(H_*,\partial H_*)\to C_*(\hat U,\partial \hat U)$. According to \cite{13} if this extends to a simple chain equivalence of free based complexes $D_*\to C_*(\hat U,\partial \hat U)$ such that $D_*$ is nonzero only in degrees 2, 3, and 4, then $\hat U$ decomposes as a union of 2-handlebodies $X\cup Y$ with $\partial \hat U\subset X$ and $X$ has spine homotopic (rel boundary) to $\partial \hat U\cup H_*$, and $Y\subset \hat U$ is a $\pi_1$ equivalence. The 2-handles in $X$ are the 2-handles needed to construct the core. It therefore remains to find a suitable chain complex.

Recall that the 1-complex $Y$ is included in $\hat N_0$, so there is a chain map 
$$C_*(M_0,Y)\to C_*(M_0, \hat N_0)\sim C_*(\hat U,\partial \hat U).$$
Next recall that  the handles $H_*$ were chosen so the inclusion  $Y\to \hat N_0\cup H_*$ is a simple homotopy equivalence. This implies the sum 
$$C_*(M_0,Y)\oplus C_2(H_*)\to  C_*(\hat U,\partial \hat U)$$
is a chain equivalence, and a minor elaboration shows it is simple. Finally since $Y\to M_0$ is surjective on $\pi_1$ it follows that $C_*(M_0,Y)$ is equivalent to a complex  concentrated in degrees 2, 3, and 4. This gives the complex needed as a template for the dual decomposition, and completes the proof.

\head 3. Properties of cores \endhead
A fair amount of structure is packed into cores to make them easier to manipulate. For example we use smooth cores even in topological s-cobordisms in order to make handlebody theory (such as it is in these dimensions) available. Explanations of the ingredients are given in 3.1. A product result in 3.2 gives the first application of the boundary spine hypothesis. This is used to give product s-cobordisms core decompositions with arbitrary cores in 3.4. Graph-sums are defined in 3.5 and used to modify core decompositons in~3.7.

\subhead 3.1 Comments on the definition\endsubhead 
The regular neighborhoods $W^5$ used in 2.1(1) are standard and determined by the orientation homomorphism $\omega\:\pi_1\to Z/2$: they are boundary connected sums of either $D^4\times S^1$ or the nonorientable analog $D^4\tilde{\times} S^1$. More generally the notation $W^n$ is used for the $n$-dimensional regular neighborhood $n\geq2$ determined by $\omega$ (or $W^n_{\omega}$, if $\omega$ is not understood from the context). Note $W^n\simeq W^2\times D^{n-2}$. 

The spine hypotheses in 2.1(3) imply in particular that each $N_i$ has the homotopy type of a 1-complex. It follows that the inclusions into $W^5$ are homotopy equivalences. Since the Whitehead group of free groups is trivial these are automatically simple. This means $W^5$ is in fact an s-cobordism (rel $\partial N$) from $N_0$ to $N_1$. 

The fundamental group of the boundary is a rough measure of the complexity of these $N$. More precisely let $W^4\subset N$ be a regular neighborhood of a 1-complex in $N$ that induces isomorphism on fundamental group. This is well-defined up to isotopy, because homotopy implies isotopy for 1-complexes in 4-manifolds. Inclusion of $\partial W^4$ into the region between $\partial N$ and $\partial W^4$ is a homotopy equivalence. An inverse gives a map of pairs $(N,\partial N)\to (W^4,\partial W^4)$. This is a homotopy equivalence $N\to W^4$ but on the boundary only a homology equivalence with $Z[\pi_1N]$ coefficients. The fundamental group of $\partial N$ is typically much larger.

Homotopy 2-deformation (or ``Andrews-Curtis equivalence'') used in 2.1(3) is the equivalence relation on 2-complexes obtained by adding or deleting cancelling pairs of cells of dimensions $(0,1)$ or $(1,2)$, and changing attaching maps of 2-cells by homotopy. The next lemma illustrates the significance of the deformation hypothesis. 

\proclaim{3.2 Lemma} Suppose $N$ is a compact 4-dimensional 2-handlebody whose spine homotopy 2-deforms to a 1-complex, and let $W^4\subset N$ be a regular neighborhood of such a 1-complex in the interior. Then $W^4\times\{0\}\subset N\times I$ extends to a diffeomorphism $W^4\times I\simeq N\times I$.  In particular $(N\times I,N\times\{0\})$ is a core s-cobordism.\endproclaim
\subsubhead Indication of proof\endsubsubhead
$(N\times I,W^4\times\{0\})$ is a relative 5-dimensional 2-handlebody with the same spine as $N$, so it homotopy 2-deforms rel $W^4\times\{0\}$ to $W^4\times\{0\}$. But in dimension 5 a homotopy 2-deformation of the spine can be realized by geometric handlebody moves. This was essentially shown by Andrews-Curtis, though their explicit result is that if the spine 2-deforms to a point then the 5-manifold is diffeomorphic to the 5-ball.
The key point is that attaching maps of 2-handles are circles in 4-manifolds, and homotopy implies isotopy for 1-manifolds in a 4-manifold. In this case the relative 2-deformation of the spine to the trivial relative 2-complex gives a handle deformation of the manifold to a collar on $W^4\times\{0\}$. This completes the proof of the Lemma.

A consequence of Lemma 3.2 is that cores can be trivialized by adding an appropriate  product s-cobordism along the boundary. This refines an idea of Matveyev \cite{10}.
\proclaim{3.3 Proposition} Suppose $(W^5,N_0,N_1)$ is a core s-cobordism, and suppose $W^4\subset N_0$ is a neighborhood of a 1-complex homotopy equivalent to $N_0$. Then attaching the product s-cobordism $(N-\text{int}\,W^4)\times I$ along the boundary s-cobordism $\partial N\times I$ gives an s-cobordism diffeomorphic to the product $(W^4\times I, W^4\times\{0\})$.
\endproclaim
In other words the product s-cobordism $W^4\times I$ has another core-product decomposition with the given core. 

\subsubhead Proof\endsubsubhead First note that adding a collar on $N_0\subset W^5$ gives again $W^5$. Specifically add $N_0\times[0,2]$ by identifying $N_0\times\{0\}$ with $N_0$. Next regard $N_0$ as the union of a smaller copy $\bar N_0$ and a collar on the boundary. Now we can identify the extended s-cobordism of the Proposition as a subset of $W^5\cup_{N_0}N_0\times[0,2]$, namely the complement of the interior of $(W^4\times [1,2])\cup (\bar N_0\times [0,1])$. Under this identification the boundary $N_0\times\{0\}\cup_{\partial N}(N_0-\text{int}\,W^4)\times \{0\}$ in the extended s-cobordism goes to the intersection of this complement with the boundary of $W^5\cup_{N_0}N_0\times[0,2]$. 

According to Lemma 3.2 the pair $\bigl((W^4\times [1,2])\cup (\bar N_0\times [0,1]),W^4\times\{2\}\bigr)$ is diffeomorphic to a collar on $W^4\times I$. But by uniqueness of collars, the complement is diffeomorphic to the original manifold. Therefore the extended s-cobordism is diffeomorphic to $W^5$, by a diffeomorphism that takes the indicated boundary to the complement of $W^4\subset\partial W^5$. But $(W^5,W^4)$ is diffeomorphic to the product $(W^4\times I,W^4\times\{0\})$. Composing these diffeomorphisms gives the diffeomorphism needed for the Proposition.

\proclaim{3.4 Corollary}Suppose $W^4\subset M^4$ is a regular neighborhood of a 1-complex and $(W^5,N_0,N_1)$ is a core s-cobordism with total space  $W^5=W^4\times I$. Then there is a decomposition of the product s-cobordism $M\times I$ with core diffeomorphic to the given one.\endproclaim
\subsubhead Proof\endsubsubhead $M\times I$ contains the product core $W^4\times I$, and according to 3.3 this decomposes as the given core and a product. 

\subhead 3.5 Graph-sums of cores\endsubhead 
This is an operation that joins two core s-cobordisms with the same fundamental group. Thinking of the boundary fundamental group as a measure of cores, we see this operation greatly increases complexity. It is used to get arbitrarily bad cores in arbitrary s-cobordisms.

Suppose $N_0\cup N_1\simeq  \partial W^5\simeq  N'_0\cup N'_1$ are two core s-cobordisms, and suppose $j\:W^3\to \partial N$, $j'\:W^3\to \partial N'$ are embeddings that induce isomorphism to the fundamental groups of $N$, $N'$. These embeddings are regular neighborhoods of graphs in a 3-manifold. They always exist but are highly nonunique because they can be knotted. Extend these embeddings to embeddings $W^3\times I\subset \partial W^5$. 

The {\it graph-sum\/} of the cores along these embeddings is $W^5\cup_{W^3\times I}W^5$ with boundary pieces $N_0\cup_{W^3}N'_0$ and $N_1\cup_{W^3}N'_1$. 

\proclaim{3.6 Lemma} The graph-sum of core s-cobordisms is a core s-cobordism.\endproclaim
 \subsubhead Proof\endsubsubhead
First consider the 5-manifolds. The embeddings $W^3\times I\subset \partial W^5$ are regular neighborhoods of graphs in a 4-manifold. Graphs are determined (for these purposes) by fundamental group, so we can suppose the graph is homotopic to a standard spine of $W^5$. But homotopic 1-complexes are isotopic in a 4-manifold, and isotopic 1-complexes have isotopic neighborhoods, so the embeddings of $W^3\times I\simeq W^4$ are isotopic to the standard embeddings $W^4\times\{0\}\subset W^4\times I\simeq W^5$. This makes it clear that the union of two copies of $W^5$ over these embeddings gives another copy of $W^5$. This is condition (1) in the definition of a core. 

To verify condition (3) we must see that $N_0\cup_{W^3}N'_0$ has a 2-handlebody structure with spine that homotopy 2-deforms to a 1-complex. Rewrite the union as 
$N_0\cup_{W^3\times\{0\}}W^3\times I\cup_{W^3\times\{1\}}N'_0$. $W^3$ has a handlebody structure with only 0- and 1-handles. Therefore the product $W^3\times I$ has a handlebody structure relative to $W^3\times \{0,1\}$ with handles the product 1- and 2-handles. Adding this to the given structures on $N_0$, $N'_0$ gives a 2-handlebody structure on the union. The spine is the union of the spines of $N_0$, $N'_0$ and (spine of $W$)$\times I$ attached by the inclusion of the ends. The hypotheses on the pieces give a 2-deformation of this to 1-complexes, joined by the product. Being a little more precise, we have 1-complexes $Y$, $Y'$ obtained from the spines of $N_0$ and $N'_0$; a 1-complex spine $Z$ of $W^3$; and homotopy equivalences $Y\leftarrow Z\to Y'$. So far we have deformed the spine of the union to the union of the mapping cylinders of the two equivalences. But the mapping cylinder of a homotopy equivalence of 1-complexes can be 2-deformed to the domain, so this union 2-deforms to $Z$. This gives the spine deformation required for the definition, and completes the proof of the Lemma.

\proclaim{3.7 Corollary} If an s-cobordism has a decomposition with core $(W^5,N_0,N_1)$ then it has a decomposition with core any graph-sum of this with any other core with total space $W^5$. \endproclaim
 \subsubhead Proof\endsubsubhead
This follows from 3.3 much as 3.4 does. Begin with the hypothesized decomposition with core $(W^5,N_0,N_1)$, and the embedding $W^3\subset \partial N$ to be used in the graph-sum. A collar on this $W^3$ in the product part of the decomposition extends this to an embedding of $W^4=W^3\times I$ in the lower end of the s-cobordism. The product structure on the product part gives a further embedding of $W^4\times I$ as a sub-s-cobordism. Now apply 3.3 to describe this as a product union an arbitrary further core. We can situate the new core so its sum embedding of $W^3\times I$ is the intersection with the original core. The result is a decomposition with core the graph-sum.

\head 4. Surgery and control\endhead
An s-cobordism gives a homotopy equivalence (rel boundary) between the ends. Standard (dimension $\geq5$) surgery and s-cobordism give a well-known and useful converse: 

\proclaim{4.1 Theorem} Suppose $N_i$ are compact 4-manifolds homotopy equivalent to 1-complexes and $f\:N_0\to N_1$ is a homotopy equivalence inducing isomorphism on the boundary. Then
\roster\item there is a topological s-cobordism inducing the equivalence;
\item if $N_i$ are smooth then there is a smooth s-cobordism if and only if the dual of the Kirby-Siebenmann invariant, in $H_1(N;Z/2Z)$, vanishes; and
\item these s-cobordisms are unique up to isomorphism. In particular they are products if and only if $f$ is homotopic rel boundary to an isomorphism.
\endroster\endproclaim

For a discussion of surgery relevant to low-dimensional topology see Freedman-Quinn \cite{4 \S11}. 
\proclaim{4.2 Corollary} A core s-cobordism has a smooth, resp. topological, product structure if and only if the induced homotopy equivalence $N_0\to N_1$ is homotopic rel boundary to a diffeomorphism, resp. homeomorphism.\endproclaim
The spine and total space hypotheses are not needed here. 
 
\subhead 4.3 Controlled s-cobordisms \endsubhead
This is a controlled version of the surgery reduction 4.1. Suppose $Y$ is a 1-complex with metric in which 1-simplices are intervals of length 1 and extended by minimum path length. The use of  standard metrics gives universal (independent of particular $Y$) $\epsilon$, $\delta$  in the statements below.

A homotopy $N\times I\to Y$ is said to have {\it radius less than\/} $\delta$ if for each $x\in N$ the distance between any two points on the image of $\{x\}\times I$ is less than $\delta$.   A $\delta$ h-cobordism over $Y$ is $(W,N_0,N_1)$ with $Y\subset N_0$ and $W\to Y$ so that 
the composition $Y\to W\to Y$ is the identity and $W$ deformation retracts to $N_0$, $N_1$ and $Y$ by deformations whose compositions into $Y$ have radius $<\delta$. 

This is a little stronger than the usual notion but equivalent in the context of core-like s-cobordisms. Combining the deformations to the $N_i$ gives a $\delta$ homotopy equivalence rel boundary $N_0\to N_1$. The result is that the equivalence determines the h-cobordism just as 4.1. 

\proclaim{4.4 Theorem} Given $\epsilon>0$ there is $\delta>0$ so that if\/ $Y\subset N_0\to N_1\to Y$ satisfies:
\roster \item $N_i$ are compact 4-manifolds 
\item $N_1\to Y$ is a $\delta$ homotopy equivalence (over $Y$), and 
\item $N_0\to N_1$ is an isomorphism on the boundary and a $\delta$ homotopy equivalence over $Y$
\endroster
then there is a topological $\epsilon$ h-cobordism inducing an $\epsilon$-equivalent equivalence, and if\/ $N_i$ are smooth there is a smooth\/ $\epsilon$ h-cobordism if and only if the Kirby-Siebenmann invariant vanishes. Further if $\rho>0$ then there is $\epsilon>0$ so that any two smooth $\epsilon$ h-cobordisms from $N_0$ to $N_1$ are $\rho$ isomorphic, and any topological\/ $\epsilon$ h-cobordism has a $\rho$ product structure.\endproclaim
All except the last conclusion use high-dimensional controlled surgery and s-cobordism, see \cite{11}. The last conclusion is 4-dimensional and follows from the controlled locally-simply connected disk embedding theorem, \cite{12, 4 \S7.2}. 

Combining 4.1 and 4.4 give a criterion for topological triviality of uncontrolled (corelike) s-cobordisms:
 
 \proclaim{4.5 Corollary} There is a $\delta>0$ so that if $f\:N_0\to N_1$ is induced by an s-cobordism and there is a $\delta$ homotopy equivalence $N_1\to Y$ with $f$ a $\delta$ homotopy equivalence rel boundary over $Y$ then the original s-cobordism is topologically a product.\endproclaim
 
\subhead 4.6 Control over a 1-complex \endsubhead
We make more explicit the control conditions in 4.5 and make a guess about the key step in obtaining them.

For convenience choose a specific model for the 1-complex $Y$: the suspension of a finite set of points. This has two vertices and each 1-cell joins these two vertices.  Let $q\:Y\to I$ be the identity on each interval.

For $\delta = 1/n$ a map $g\:N\to Y$ is a $\delta$ homotopy equivalence essentially if the following holds: divide $I$ into $n$ equal subintervals. Then the $pg$ inverse image of each  subinterval $[\frac{k}{n},\frac{k+1}{n}]$ is contractible in the inverse image of the larger interval $[\frac{k-1}{n},\frac{k+2}{n}]$. 

We speculate that the key problem is getting started. 
\subhead 4.7 Question \endsubhead Suppose $N$ is a 4-dimensional 2-handlebody with spine 2-equivalent to a 1-complex $Y$ as above. Is there a homotopy equivalence $N\to Y$ so the inverse image of $[0,\frac12]$ is contractible in the inverse of $[0,\frac34]$?

There is a spiritual similarity of this to Freedman's Poincar\'e transversality formulation of the surgery problem, see the last section of \cite4. The handlebody and spine conditions are satisfied for cores so are not serious restrictions. Recall that 2-handlebodies have Kirby calculus descriptions as link diagrams, \cite{7}. The spine condition may help in manipulating these diagrams, but see \cite{6} for limitations. The induced decomposition of the 3-manifold $\partial N$ is likely to be important, and is the most likely source for obstructions if the answer to the question is ``no''.

\subhead 4.8 Question \endsubhead Suppose $N_0\to N_1$ is a homotopy equivalence of 4-dimensional 2-handlebodies with spines 2-equivalent to $Y$, and is a diffeomorphism on the boundary. Suppose $N_1\to Y$ satisfies the conclusion of Question 4.7. Are these maps homotopic to $\delta$ equivalences as required for Corollary 4.5?

The hope is that once control of $N_1\to Y$ gets started, in the sense of 4.7, then it can be improved using dual-decomposition and other techniques. Further $N_0$ and $N_1$ have the same boundary, so good control on $N_1\to Y$ gives a decomposition of $\partial N_0$ and a homotopy template for extending it to the interior. Again the hope is that existing techniques will be sufficient to promote this to $\delta$ equivalence. 

Note that the handlebody decompositions on $N_0$ and $N_1$ give two link descriptions of the boundary. Part of the question can be thought of as a weak homotopy relation between the two link descriptions, without stabilization.

\enddocument

\bye